\documentclass[conference]{IEEEtran}
\IEEEoverridecommandlockouts
\usepackage{amsmath,amssymb,amsfonts}
\usepackage{algorithmic}
\usepackage{graphicx}
\usepackage{textcomp}
\usepackage{xcolor}
\usepackage{booktabs}
\usepackage{cleveref}
\usepackage[style=ieee,backend=biber, url=false, doi=false, isbn=false]{biblatex}
\Crefname{figure}{Fig.}{Fig.}

\begin{document}

\title{Improving Numerical Stability and Accuracy in Partitioned Methods with Algebraic Prediction\\

\thanks{This work is sponsored by the National Science Foundation, Award ECCS-2339148.}
}

\author{\IEEEauthorblockN{Ahmad Ali}
\IEEEauthorblockA{\textit{Department of ECE} \\
\textit{North Carolina State University}\\
aali27@ncsu.edu}
\and
\IEEEauthorblockN{Haya Monawwar}
\IEEEauthorblockA{School of ECE\\
\textit{Oklahoma State University}\\
haya.monawwar@okstate.edu}
\and
\IEEEauthorblockN{Hantao Cui}
\IEEEauthorblockA{Department of ECE\\
\textit{North Carolina State University}\\
hcui9@ncsu.edu}
}

\maketitle

\begin{abstract}
The partitioned approach for the numerical integration of power system differential algebraic equations faces inherent numerical stability challenges due to delays between the computation of state and algebraic variables.
Such delays can compromise solution accuracy and computational efficiency, particularly in large-scale system simulations.
We present an $O(h^2)$-accurate prediction scheme for algebraic variables based on forward and backward difference formulas, applied before the correction step of numerical integration.
The scheme improves the numerical stability of the partitioned approach while maintaining computational efficiency.
Through numerical simulations on a lightly damped single machine infinite bus system and a large-scale 140-bus network, we demonstrate that the proposed method, when combined with variable time-stepping, significantly enhances the numerical stability, solution accuracy, and computational performance of the simulation. 
Results show reduced step rejections, fewer nonlinear solver iterations, and improved accuracy compared to conventional approaches, making the method particularly valuable for large-scale power system dynamic simulations.
\end{abstract}

\begin{IEEEkeywords}
Transient simulation, partitioned approach, numerical stability, differential algebraic equations
\end{IEEEkeywords}

\section{Introduction}
Transient simulation is essential for studying the dynamic behavior of power systems under different operating conditions.
Due to the critical nature of these studies, the simulation programs are expected to yield reliable results while being computationally efficient. 
At the heart of the dynamic simulation lies the numerical integration of differential-algebraic equations (DAEs), given by;
\begin{align}\label{eq:dae}
\dot{x} = f(x, y) \\
0 = g(x,y)
\end{align}
where $x \in \mathbb{R}^{m}$ is the vector of state variables, $y \in \mathbb{R}^{n}$ is the vector of algebraic variables, $f: \mathbb{R}^{m} \times \mathbb{R}^{n} \mapsto \mathbb{R}^{m}$ represents the differential equations, and
$g: \mathbb{R}^{m} \times \mathbb{R}^{n} \mapsto \mathbb{R}^{n}$ represents algebraic equations. The dependence of $f$ and $g$ on time and other simulation parameters is implied.

Numerically, the DAEs can be solved using either a simultaneous or a partitioned approach.
In the simultaneous approach, the differential and algebraic equations are solved together.
However, in the partitioned approach, these equations are solved separately in an alternating manner \cite{milano2010power}.

The simultaneous approach offers better numerical accuracy but can be computationally more expensive due to the need to rebuild and factorize a large Jacobian matrix at each step.
The partitioned approach is the \textit{de facto} standard in commercial tools.
However, it faces significant challenges due to the inherent delay between state and algebraic variables.
Specifically, when computing state variables at time $t_n$, the algebraic variables are unknown and must be approximated using values from the previous time step. 
This delay becomes particularly problematic because partitioned methods typically employ explicit numerical integration schemes, which are susceptible to numerical instability due to error accumulation.

Current implementations address this challenge by comparing the values of the algebraic variables after each time step \cite{chow2020power}.
If the updated algebraic variables are sufficiently close to the previous ones, the step is accepted.
Otherwise, the step is rejected, and additional corrector iterations are performed, accompanied by re-solving the algebraic equations.
This approach introduces two significant issues. First, even when a step is accepted, residual errors can accumulate and impact the stability of explicit methods.
Second, when steps are rejected, an algebraic loop is introduced, requiring additional Newton iterations or correction steps, or a combination of both \cite{tzounas2023unified}.
This process can be computationally expensive and can diminish the benefits of using a partitioned approach.
Moreover, the accuracy of the simulation is directly dependent on the tolerance factor chosen for the algebraic variables, which is often selected heuristically.

Previous research has identified the impact of these delays on the numerical accuracy and stability of simulation \cite{milano2016delay, tzounas2022small}.
Recently, a matrix pencil-based method was proposed \cite{tzounas2023unified} to study the stability and accuracy of the partitioned approach using the generalized eigenvalues of the system.
The authors established bounds on the step size to ensure numerical stability.
However, these methods have limitations: they rely on small signal stability analysis of the system around an equilibrium point. 
The qualitative behavior of the system, however, can not be reliably extended to large disturbances or to cases where the disturbance introduces topological changes.
Additionally, using a constant step size based on dominant eigenvalues is inefficient, since the fast dynamics settle quickly and are overtaken by slower dynamics.

The fundamental research questions that remain are:
\begin{itemize}
    \item Can a qualitative analysis be performed to identify the factors affecting the accuracy of the corrector step?
    \item Is controlling the local error alone at each integration step sufficient to ensure stability?
    \item How can the problem of ``delay'' be effectively addressed in the partitioned approach?
\end{itemize}
In this work, we perform an error analysis of the corrector step, deriving the relationship between the error term and simulation parameters, such as step size, predictor accuracy, and errors in algebraic variables.
Next, a prediction scheme is proposed to estimate the algebraic variables for the correction step.
This scheme is based on the forward and backward difference formulas derived from Taylor series expansions, ensuring simplicity in implementation and computational efficiency.
Yet, it is general such that the information regarding the numerical solver or the system model is not required \textit{a priori}.
Through numerical simulations on a poorly damped single machine infinite bus system and a large 140-bus NPCC system, it is shown that when combined with error control and variable time stepping schemes, the proposed framework ensures numerical stability while maintaining user-defined accuracy and computational efficiency.

The remainder of this paper is organized as follows:
Section II provides an overview of the Predictor-Corrector numerical integration method.
In section III, an analytical expression for the error in the corrector step is derived and the proposed algebraic approximation scheme is presented.
Section IV presents the case studies.
In section V, simulation results are discussed.
Finally, section VI concludes the paper. 
\section{Methodology}
\subsection{Predictor Corrector Scheme}
The Predictor-Corrector (PC) method is a type of linear multi-step method that combines the advantages of both explicit and implicit methods.
While explicit methods are simple to compute, they suffer from poor stability.
Implicit methods, on the other hand, are more stable but require iterative solutions, which can be computationally expensive.

The PC method addresses this by first using an explicit method to estimate the solution and then 
an implicit method is used to refine this estimate.
This approach avoids the need for repeated iterations of the implicit equation while improving stability, making it a popular choice for power system simulation tools.

\subsubsection{Predictor Step}
In this step, an explicit method is used to estimate the state variables at $t_{n}$.
An explicit method approximates the successive values of $x_{n+1}$ in terms of the previously computed values \cite{ascher2011first}.
Such methods can be represented in a general form as:
\begin{equation}
    x_{n+1} = \text{H}(f,t_{n},\dots ,t_{n+1-m}, x_n,\dots , x_{n+1-m})
\end{equation}
where $\text{H}$ represents the formula for the explicit method and $m>0$ is the number of previous steps considered in the method.

Since all the quantities required to carry out this computation are available, the prediction can be carried out in a straightforward way.
Depending on the order of accuracy, $q$, of the explicit method, an $O(h^q)$ accurate estimate of $x_{n+1}$ is obtained, which is represented as $\Tilde{x}_{n+1}$.

\subsubsection{Corrector Step}
In the corrector step, an implicit method is used to approximate $x_{n+1}$ based on $\Tilde{x}_{n+1}$.
An implicit method defines $x_{n+1}$ as the solution of an equation that depends implicitly on $x_{n+1}$:
\begin{equation}\label{eq:implicit_formula}
    0 = \text{G}(f,t_{n+1},\dots ,t_{n+1-m}, x_{n+1},\dots , x_{n+1-m})
\end{equation}
where $\text{G}$ represents the implicit formula and $m>0$ represents the number of previous steps considered in the method.
\eqref{eq:implicit_formula} can be written in the form of a fixed-point equation as;
\begin{align}\label{eq:fixed_point_1}
    x_{n+1} = \Phi(x_{n+1})
\end{align}
In the PC scheme,
only a fixed number of iterations of the corrector step are performed, using $\Tilde{x}_{n+1}$ as the initial estimate.

As shown in \eqref{eq:dae}, the differential equations depend on the algebraic variables as well.
Thus, \eqref{eq:fixed_point_1} can be modified as:
\begin{align}\label{eq:fixed_point_1}
    x_{n+1} = \Phi(x_{n+1}, y_{n+1})
\end{align}
where $y_{n+1}$ is a parameter in this formula.
Since, 
is unknown at this time, as an estimate is used to solve \eqref{eq:fixed_point_1}.

In power system simulation software, generally, the value of the algebraic variables at the previous time step, $y_{n}$, is used as an estimate for $y_{n+1}$.
However, such an approach introduces a ``delay'' between the state variables and the algebraic variables.
The errors associated with the delay are explored in the following section.
\subsection{Computation of Algebraic Variables}
After the correction step, the algebraic equations are solved for $y_{n+1}$ using the corrected values of state variables, $x_{n+1}$.
\begin{align}
    0 = g(x_{n+1}, y_{n+1})
\end{align}

In this work, the forward Euler method is used as the predictor while the implicit trapezoidal method is used as the corrector.

\section{Error Analysis and Proposed Scheme}
\subsection{Impact of Using Old Values of the Algebraic Variables on Accuracy of Numerical Integration}
To investigate the impact of using $y_n$ as an approximation for $y_{n+1}$ in the correction step, a mathematical framework is derived.
Without loss of generality, assume that the backward Euler (BE) method is used in the correction step.
The analysis that follows is valid for any implicit numerical integration method.
Using \eqref{eq:dae},
the difference equation obtained using BE is given as;
\begin{equation}\label{eq:be_formula}
    x_{n+1} = x_{n} + hf(x_{n+1}, y_{n+1})
\end{equation}
where the subscript $n$ represents the integration step, $h$ is the step size and $x_n$ and $y_n$ are the values of state and algebraic variables at step $n$, respectively.
Assuming that the corrector iteration is solved until convergence.
Then;
\begin{align}\label{eq:be_formula2}
    x_{n+1}^{1} = x_{n} + hf(x_{n+1}^{1}, y_{n+1}^{0})
\end{align}
where $x_{n+1}^{1}$ represents the solution obtained using the corrector iteration
and $y_{n+1}^{0}$ is the approximation of $y_{n+1}$.

Let $x_{n+1}^{*}$ and $y_{n+1}^{*}$ be the true solutions of the system at $t_{n+1}$.
Then, the error in $x_{n+1}^{1}$ is:
\begin{align} \label{eq:e_x_initial}
    e_x = x_{n+1}^{1} - x_{n+1}^{*}
\end{align}
Using \eqref{eq:be_formula} and \eqref{eq:be_formula2}, we get:
\begin{align}
    e_x &= [x_{n} + hf(x_{n+1}^{1}, y_{n+1}^{0})] - [x_{n} + hf(x_{n+1}^{*}, y_{n+1}^{*})] \\
    &= h [f(x_{n+1}^{1}, y_{n+1}^{0}) - f(x_{n+1}^{*}, y_{n+1}^{*})] 
\end{align}
To find this difference, we can use the mean value theorem (MVT) for multi-variable functions.
Application of MVT results in:
\begin{equation}\begin{split}\label{eq:MVT}
    f(x_{n+1}^{1}, y_{n+1}^{0}) - f(x_{n+1}^{*}, y_{n+1}^{*}) &= \\ f_x (x_{n+1}^{1} - x_{n+1}^{*}) + f_y (y_{n+1}^{0} - y_{n+1}^{*})
\end{split}\end{equation}
where $f_x$ and $f_y$ are the partial derivatives of $f(x,y)$ at some point between $(x_{n+1}^{1}, y_{n+1}^{0})$ and $(x_{n+1}^{*}, y_{n+1}^{*})$, respectively.
Let $e_y = y_{n+1}^{0} - y_{n+1}^{*}$, then;
\begin{align}\label{eq:e_x_initial}
    e_x = h(f_x e_x + f_y e_y)
\end{align}
Rearranging;
\begin{align}\label{eq:e_x_final}
    e_x = \frac{hf_y}{1 - hf_x} e_y
\end{align}
\eqref{eq:e_x_final} reveals that the error in the estimate of $x_{n+1}$ is directly proportional to the error in $y$ scaled by the step size $h$.
When $e_y$ is uncontrollable, then the only option is to use a small $h$, which is computationally sub-optimal.

Moreover, when the step size is not small enough, the errors continue to increase in an unbounded manner,
since the error in $x_{n+1}$ introduces error in $y_{n+1}$, creating a feedback loop.
To avoid this problem, a common strategy is to compare $y_{n+1}$ with $y_n$ once the algebraic equations have been solved for time $t_{n+1}$.
If this difference exceeds a set threshold, $\epsilon$, the corrector iterations are repeated followed by the solution of the algebraic equations.
This process is continued until some accuracy criteria is met.
However, this introduces undesired computational expense.
Moreover, $\epsilon$ is often set heuristically, which can impact the numerical accuracy and computational stability of the simulation.

When $y_{n+1}^{0} = y_n$ is used in \eqref{eq:be_formula2}, then;
\begin{align}\label{eq:e_y_initial}
    e_y = y_n - y_{n+1}^{*}
\end{align}
Assume that $y_n$ is the exact solution at $t_n$.
Further, if we assume $y$ to be sufficiently smooth in a small neighborhood around $y_n$, 
then using Taylor series expansion;
\begin{align}
    y_{n+1}^{*} = y_{n} + h y_{n}^{'} + O(h^2)
\end{align}
Then \eqref{eq:e_y_initial} becomes;
\begin{align}
    e_y &= y_n - y_{n+1}^{*} = -hy_{n}^{'} - O(h^2)
\end{align}
Using the leading term in $e_y$, we get;
\begin{align}
    e_y = O(h)
\end{align}
Thus, the error introduced by using $y_{n}$ as an estimate for $y_{n+1}$ in the correction step decreases linearly with time.

In the analysis above, we assumed that the corrector iteration is solved until convergence.
However, if only a single iteration of the corrector is applied, then we replace $x_{n+1}^{1}$ in the right-hand side of \eqref{eq:MVT} with the solution of the predictor step, $\Tilde{x}_{n+1}$. 
Then after solving for $e_x$, we get:
\begin{align}
    e_x \approx h(f_x O(h^q) + f_y e_y)
\end{align}
where $q$ is the order of accuracy of the predictor method.
This reveals that when the correction iteration is not solved until convergence, then the error in $x_{n+1}$ is dependent on the order of accuracy of the predictor method.

\subsection{Proposed Algebraic Prediction Scheme}
The problem of unavailability of $y_{n+1}$ can be addressed by making an approximation of $y_{n+1}$ that reduces the contribution of $e_y$ in \eqref{eq:e_x_final}.
Note that at the start of the integration step at $t_{n+1}$, we have the knowledge of the previous values of $y$, i.e. $y_{n}$, $y_{n-1}$ etc.
The goal is to use this knowledge to obtain a cheap yet better approximation of $y_{n+1}$, without solving the non-linear algebraic equations.
The idea is based on the Taylor series expansion of $y_{n+1}$ around $y_n$.

Given that the underlying function, $y$, is sufficiently smooth in a small neighborhood around $y_n$, then the Taylor series expansion is given as:
\begin{equation}\label{eq:taylor_series}
    y_{n+1} = y_n + hy'_{n} + O(h^2)
\end{equation}
To find $y'_{n}$, we use the two-point backward difference formula; 
\begin{equation}\label{eq:taylor_diff}
        y'_{n} \approx \frac{y_n - y_{n-1}}{h}
\end{equation}
Substitute \eqref{eq:taylor_diff} in \eqref{eq:taylor_series} and neglecting the higher order terms, we get an estimate for $y_{n+1}$, as:
\begin{align}\label{eq:y_approx}
    \Tilde{y}_{n+1} \approx y_n + h_{n+1} ( \frac{y_n - y_{n-1}}{h_{n}})
\end{align}
where, $h_{n+1}$ and $h_n$ are the current and previous time steps, respectively.
The obtained approximation of $y_{n+1}$ is $O(h^2)$-accurate.

The above analysis is based on the assumption that the underlying function is sufficiently smooth in a small neighborhood.
However, this assumption is violated at the switching instants, during which the algebraic variables undergo large, often discontinuous changes.
At such instances, this approximation can be dropped and the standard approximation scheme can be employed.
However, after the switching instants, the proposed scheme can be effectively used.

To estimate $e_y$ when the proposed scheme is used, we use \eqref{eq:taylor_series} and \eqref{eq:y_approx} to get;
\begin{align}
    e_y = \Tilde{y}_{n+1} - y_{n+1} = O(h^2)
\end{align}
Thus, by using the proposed strategy to approximate $y_{n+1}$ in the corrector step, the numerical accuracy of the corrector step can be enhanced without the need to reduce the step size.
The proposed approach is simple to implement and is computationally inexpensive, as it only involves vector arithmetic operations, which are very fast.  

\section{Case Studies}
Two case studies are conducted to investigate the impact of the proposed scheme on the numerical accuracy and stability of transient simulation using PC approach.
The first case study is based on a poorly damped single-machine infinite bus (SMIB) system taken from \cite{chow2020power}.
The second case study is based on a 140-bus NPCC system with turbine governors and generator exciters, representing a fairly damped system.

The PC solver is implemented in the open-source power system simulation tool ANDES \cite{cui2020hybrid}, which natively employs a simultaneous approach using the implicit trapezoidal method (ITM) for numerical integration.
The ITM-based results serve as benchmarks.
A proportional-integral (PI)-based adaptive time stepping controller is used to adjust the step size based on the estimate of local error.

The proposed prediction-based PC scheme is compared with the standard approach that estimates $y_{n+1}$ using $y_{n}$ in the correction stage.
In the standard approach, large differences in algebraic variables lead to a reduction in the step size to maintain small differences between successive steps.

Computational performance is assessed by analyzing the simulation parameters such as integration step sizes following disturbances, the total number of non-linear solver calls, and the number of accepted and rejected steps.

\subsection{Case Study 1: SMIB System}
The system consists of a synchronous generation at bus 1, feeding the slack bus at bus 2.
Bus 3 is modeled as a PQ bus.
At $t=0.5$ seconds, a three-phase line-to-ground fault is applied at bus 3 and cleared after $0.1$ seconds.

\begin{figure}
    \centering
    \includegraphics[scale=0.45]{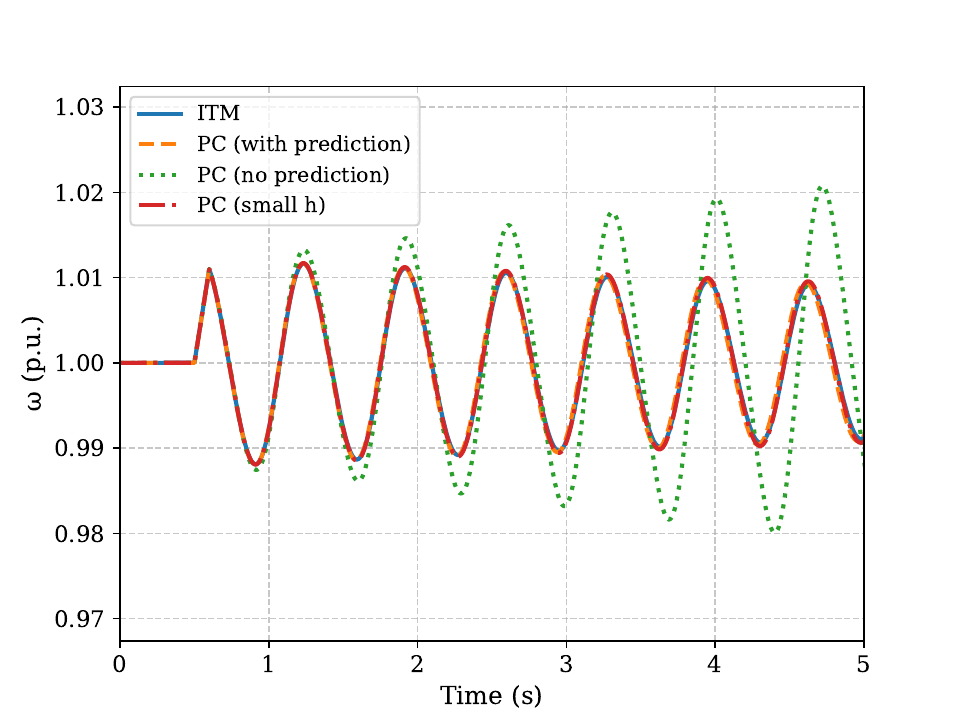}
    \caption{Comparison of Generator 1 speed for the SMIB system.}
    \label{fig:smib_PC_unstable}
\end{figure}

\Cref{fig:smib_PC_unstable} shows the trajectory of the generator speed, ($\omega$), obtained using the ITM and the PC approaches.
As observed, the results from the ITM and the PC approach with the proposed prediction scheme closely agree. 
In contrast, the standard PC approach using $y_{n}$ as an estimate for $y_{n+1}$ exhibits diverging behavior, indicative of numerical instability.

To stabilize the standard PC approach, a tighter bound on the allowable difference between $y_{n}$ and $y_{n+1}$ was introduced.
In this case, the PI-based step controller reduced the step size to keep differences within acceptable limits.
While it improved the stability, the computational burden was increased due to the requirement of very small step sizes.

\begin{table}[]
\caption{Computational performance metrics of the solvers for SMIB}
\label{tab:smib_stats}
\begin{tabular}{@{}cccc@{}}
\toprule
Solver                 & Nonlinear Calls & Accepted Steps & Rej. Steps \\ \midrule
ITM                    & 2,318             & 759            & 59             \\
PC (No Prediction)    & 2,904             & 900            & 89             \\
PC (With Prediction)      & 3,247             & 763            & 90             \\
PC (Step Adjustment) & 30,502            & 15,030          & 254            \\ \bottomrule
\end{tabular}
\end{table}

\Cref{tab:smib_stats} presents the computational metrics for the solvers over a 10-second simulation.
While the standard PC solver with tighter bounds achieves accuracy comparable to the ITM and the PC with prediction, it incurs a significantly higher number of non-linear solver calls.
The PC approach with prediction attains similar accuracy and stability with slightly more nonlinear solver calls than the ITM, indicating the use of a smaller time step than the ITM.
\begin{figure}
    \centering
    \includegraphics[scale=0.37]{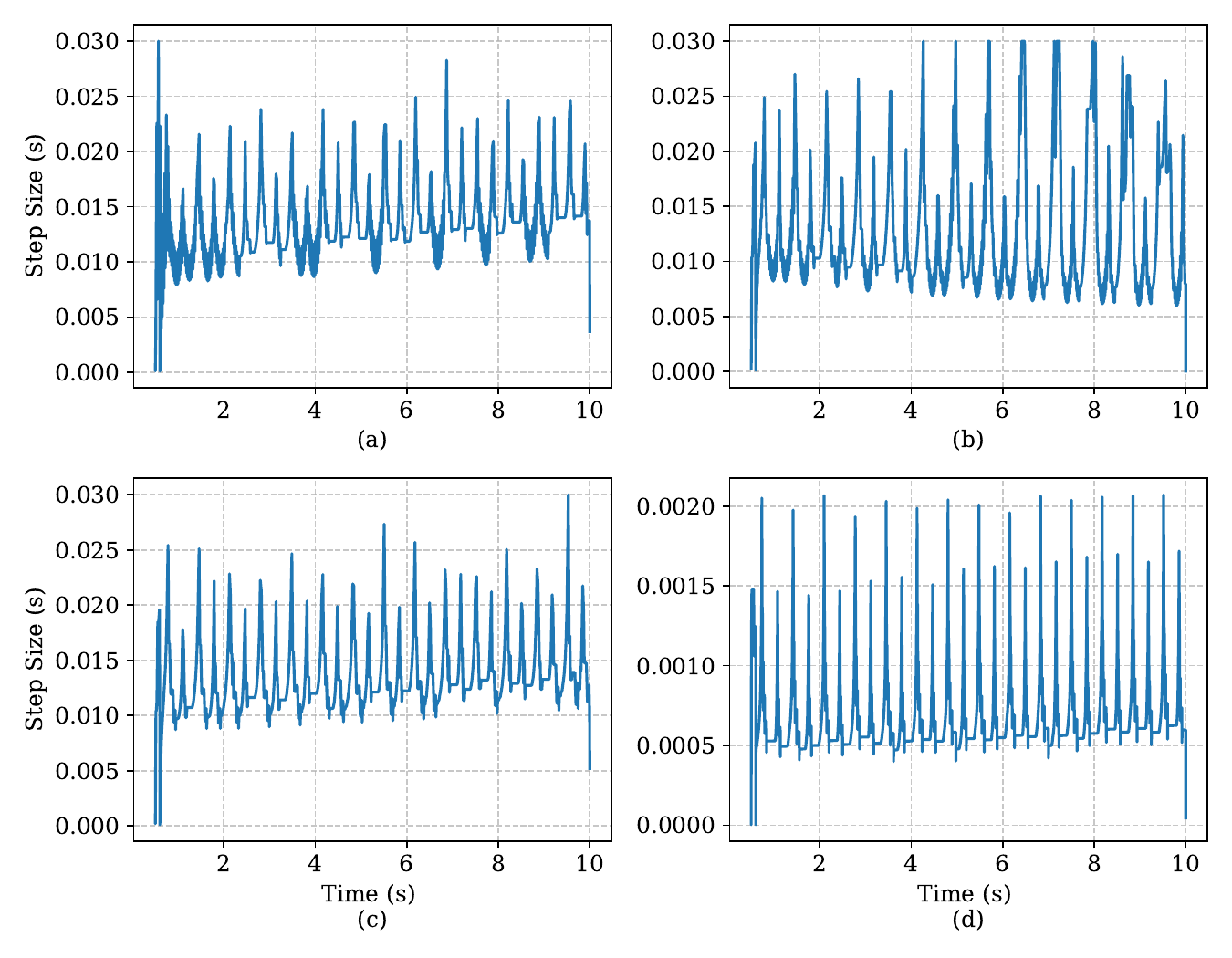}
    \caption{Step size used by the integration scheme. (a) ITM (b) PC approach with no prediction (unstable case) (c) PC approach with prediction, (d) PC with no prediction (stable case)}
    \label{fig:smib_dt}
\end{figure}

\Cref{fig:smib_dt} illustrates the step sizes used in the numerical integration for each scheme following the disturbance.
As observed, the step sizes employed by the prediction-based PC approach are comparable to those of ITM.

These results demonstrate that controlling local error alone cannot ensure numerical stability.
The standard PC converges with tighter error bounds because the tighter bounds force the solver to take very small steps, as shown in \cref{fig:smib_dt}, which contributed to the numerical stability of the solver.
As indicated from our analysis of $e_x$, we observe that to maintain the same level of accuracy, the standard PC approach had to use a step size at least an order of magnitude smaller than that used by the PC approach with prediction.
\subsection{Case Study 2: NPCC System}
In this case study, a line trip fault is simulated.
Line 2, connecting buses 1 and 4, is tripped at $t=0.5$ seconds and reconnected at $t=0.6$ seconds.
No numerical instability was observed during this simulation.
\begin{figure}
    \centering
    \includegraphics[scale=0.43]{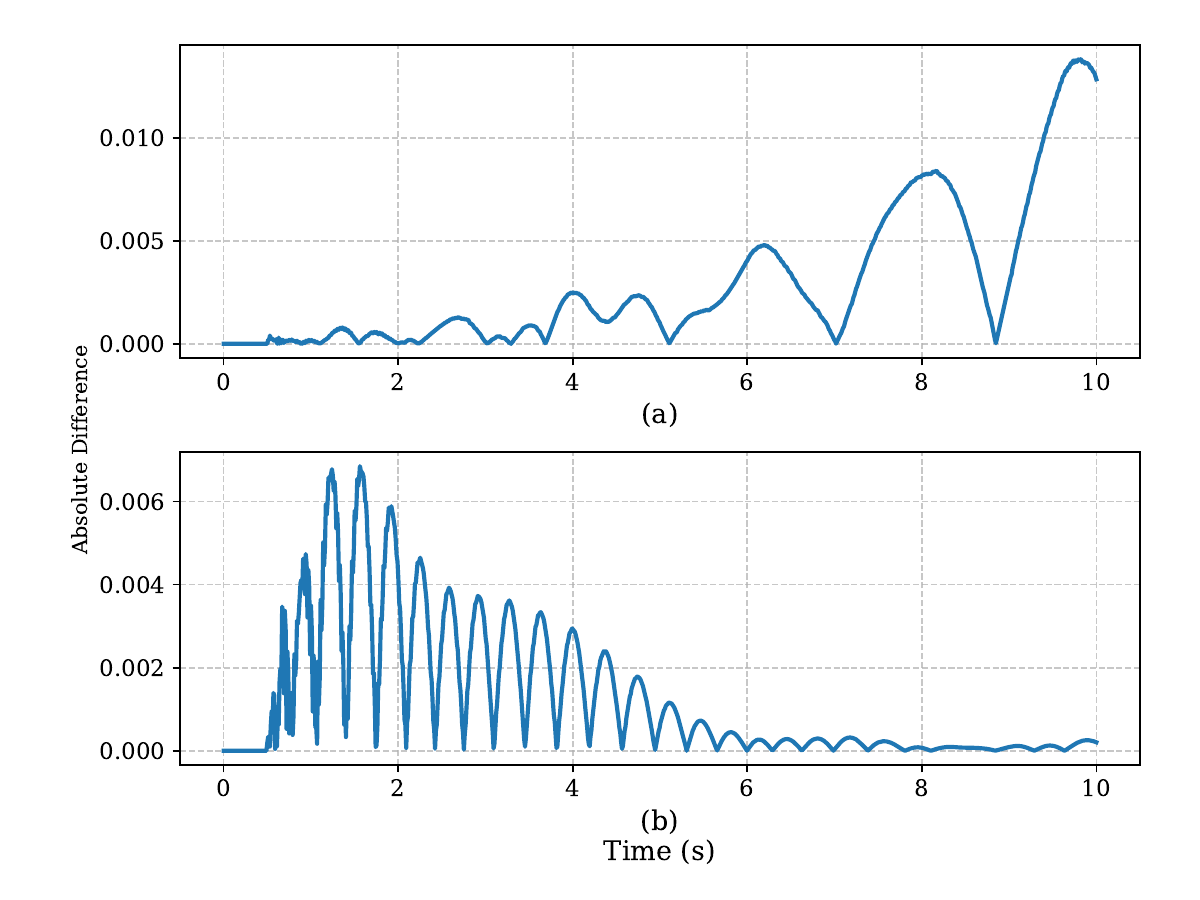}
    \caption{Absolute difference in the simulation results. (a) PC approach with no prediction, (b) PC approach with prediction}
    \label{fig:npcc_diff}
\end{figure}

Simulations were conducted using both the standard PC scheme and the prediction-based PC scheme, with results compared to those obtained using the ITM solver.
For each variable, the absolute difference between the results obtained from each method and the ITM was computed.
The variable with the maximum difference, determined using $L_2$ norm, was identified for each method.
The maximum differences for each method are illustrated in \Cref{fig:npcc_diff}.

\begin{figure}
    \centering
    \includegraphics[scale=0.45]{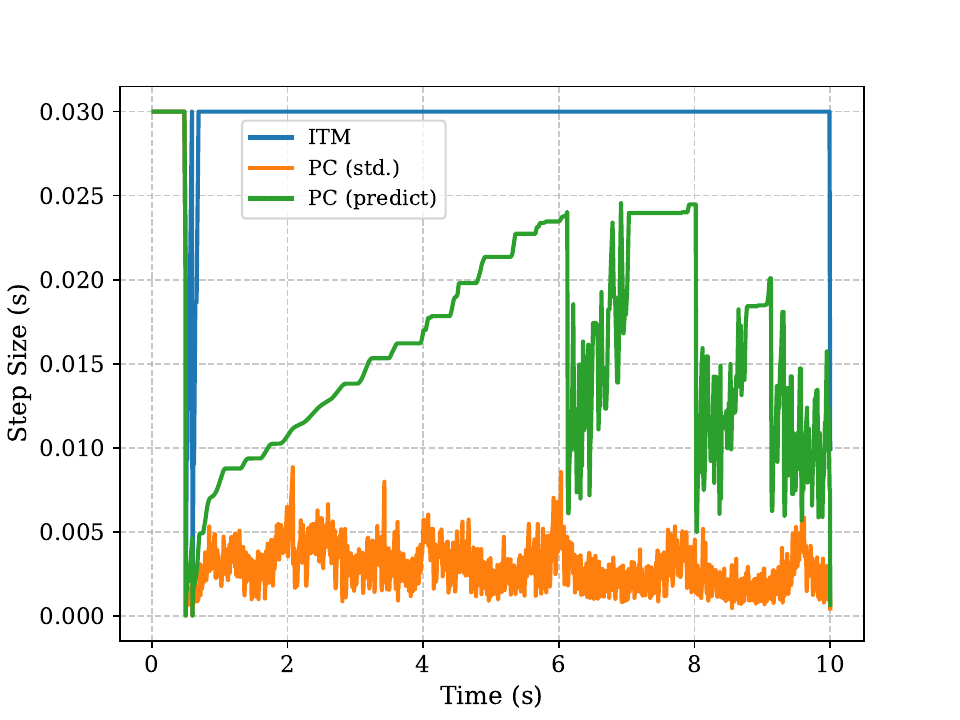}
    \caption{Step size used in numerical integration for the simulation of NPCC.}
    \label{fig:npcc_dt}
\end{figure}
\Cref{fig:npcc_dt} compares the step sizes used by the solvers.
We observe that the standard PC implementation requires a step size at least an order of magnitude smaller than that used by the prediction-based PC implementation

\begin{table}
\caption{Computational Performance Metrics of the Solvers for NPCC}
\label{tab:npcc_stats}
\begin{tabular}{@{}cccc@{}}
\toprule
Solver              & Nonlinear Calls & Accepted Steps & Rej. Steps \\ \midrule
ITM                 & 692              & 350            & 0              \\
PC (No Prediction) & 12,089            & 3,880           & 848            \\
PC (With Prediction)   & 2,566             & 804            & 66             \\ \bottomrule
\end{tabular}
\end{table}
\Cref{tab:npcc_stats} summarizes the computational performance metrics of the solvers.
As observed, the prediction-based PC scheme outperforms the standard PC scheme.
However, the simultaneous ITM method surpasses both partitioned approaches, achieving the fewest number of calls to the nonlinear solver and zero rejected steps.

\section{Conclusion}
A prediction scheme for algebraic variables is proposed to improve the numerical accuracy and stability of the predictor corrector-based partitioned approach for solving power system DAEs.
The scheme employs forward and backward difference formulas to obtain an $O(h^2)$-accurate approximation of the algebraic variables used in the correction stage of the PC scheme.
Key conclusions of this work are:
\begin{itemize}
    \item The error incurred in the correction step is directly proportional to the error in the algebraic variables, the order of accuracy of the prediction method, and the simulation step size.
    \item Numerical results demonstrate that controlling the local error alone does not ensure numerical stability.
    \item Compared to the standard PC approach, the proposed scheme reduces the error introduced by the algebraic variables in the differential equations from $O(h)$ to $O(h^2)$.
\end{itemize}

\printbibliography
\end{document}